\documentclass[11pt,reqno]{amsart}
\usepackage[T1]{fontenc}
\usepackage[utf8]{inputenc}
\usepackage{lmodern}
\usepackage{amsmath,amssymb,amsthm,mathtools}
\usepackage[a4paper,textwidth=150mm,textheight=230mm,centering]{geometry}
\usepackage{microtype}
\usepackage[hidelinks]{hyperref}
\hypersetup{pdftitle={Tensor Products and the Isomorphic Theory of Daugavet Spaces},
 pdfsubject={Injective and projective tensor products of inductively constructed Banach spaces}}
\numberwithin{equation}{section}
\newtheorem{theorem}{Theorem}[section]
\newtheorem{lemma}[theorem]{Lemma}
\newtheorem{proposition}[theorem]{Proposition}
\newtheorem{corollary}[theorem]{Corollary}
\theoremstyle{remark}

\newcommand{\R}{\mathbb R}

\newcommand{\eps}{\varepsilon}
\newcommand{\id}{\operatorname{Id}}
\newcommand{\tr}{\operatorname{tr}}

\newcommand{\av}{\mathbb E}
\newcommand{\te}{\widehat\otimes_{\eps}}
\newcommand{\tp}{\widehat\otimes_{\pi}}
\newcommand{\norm}[1]{\left\lVert #1\right\rVert}
\newcommand{\abs}[1]{\left\lvert #1\right\rvert}
\newcommand{\ip}[2]{\left\langle #1,#2\right\rangle}
\newcommand{\Slice}{\operatorname{Slice}}
\allowdisplaybreaks[2]
\title{Tensor Products and the Isomorphic Theory of Daugavet Spaces}
\author{Samir Hamad}
\address{Mathematisches Seminar, Christian-Albrechts-Universit\"at zu Kiel, Kiel, Germany}
\email{hamad@math.uni-kiel.de}
\date{October 2026}
\subjclass[2020]{Primary 46B04, 46B28; Secondary 46B20}
\keywords{Daugavet property, injective tensor product, projective tensor product,
 strongly exposed point, inductive construction, Lipschitz function spaces, Renormings}
\begin{document}
\begin{abstract}
We construct two separable Banach spaces with the Daugavet property whose
injective and projective tensor squares, respectively, fail the Daugavet
property. Both examples arise from one sequence of finite-dimensional
extensions. The first space is the completion of the increasing union of
the extensions, whereas the second is the closed span of their coordinate
functionals. We will also investigate the operator Daugavet properties of these spaces later in the paper, also highlighting natural examples of metric spaces \(M\) for which \(\operatorname{Lip}_0(M)\) has the Daugavet property but fails to have the weak operator Daugavet property. Independently of this, we will also prove that every Banach space admits a linear isometric embedding into
    a Banach space which is not isomorphic to any Banach space with the
    Daugavet property. 
\end{abstract}
\maketitle

\section{Introduction}
Throughout the paper all spaces, scalars and operators are real. We say that a normed space $X$ has the \emph{Daugavet property}, shortly $X \in \mathrm{DPr}$, if
\[
\|\mathrm{Id} + T\| = 1 + \|T\|,
\]
for every rank-one operator $T \colon X \to X$. The following characterization was obtained in \cite[Lemma~2.2]{KSSW}.

\begin{theorem}
Let $X$ be a normed space. The following assertions are equivalent:
\begin{enumerate}
\item $X$ has the Daugavet property.
\item For every $y \in S_X$, $x^* \in S_{X^*}$, and every $\varepsilon > 0$, there exists $x \in S_X$ such that
\[
x^*(x) \ge 1 - \varepsilon
\quad \text{and} \quad
\|x + y\| \ge 2 - \varepsilon.
\]
\item For every $\varepsilon > 0$ and every $y \in S_X$, the closed convex hull of
\[
\{\, u \in (1+\varepsilon) B_X \mid \|y + u\| \ge 2 - \varepsilon \,\}
\]
contains $S_X$.
\end{enumerate}
\end{theorem}

Whether the Daugavet property passes from two factors to their injective
or projective tensor product is a longstanding question; see
\cite[Section~6]{Werner}.
The version with only one Daugavet factor was answered negatively by
Kadets, Kalton and Werner in the complex setting \cite[Theorem~4.2 and Corollary~4.3]{KKW}.
Several positive results are known under additional assumptions.
Rueda Zoca, Tradacete and Villanueva proved that the injective tensor
product of two atomless $L_1$-spaces has the Daugavet property, and that
the projective tensor product of two $L_1$-preduals with the Daugavet
property also has this property \cite{RTV}. Mart\'in and Rueda Zoca
introduced the weak operator Daugavet property and proved its stability
under projective tensor products \cite[Theorem~5.4]{MRZ}. For the injective
tensor product, Rueda Zoca proved that it suffices to assume the weak
operator Daugavet property of both dual spaces \cite[Theorem~3.5]{RZ}. We show that the Daugavet property need not be preserved under taking injective or projective tensor squares. In both cases, our approach is based on the standard technique of constructing a sequence of finite-dimensional extensions. 

The last part of this paper will consider a short result regarding the isomorphic theory of Daugavet spaces. It is a known result \cite[Corollary~8.4.12]{KadetsMartinRuedaZocaWerner} that every separable
Banach space containing an isomorphic copy of
$C[0,1]$ admits an equivalent norm under which
it has the Daugavet property.
Thus, although the Daugavet property is not
invariant under equivalent renormings, the
presence of an isomorphic copy of $C[0,1]$
provides a sufficient condition for a separable
Banach space to be renormed to have the
Daugavet property. We will show that this is not true for nonseparable Banach spaces.

We recall the tensor notation. For $u=\sum_{i=1}^q x_i\otimes y_i\in
E\otimes F$, put
\[
 \norm u_\eps=
 \sup_{x^*\in B_{E^*},\,y^*\in B_{F^*}}
 \abs{\sum_{i=1}^q x^*(x_i)y^*(y_i)},
 \qquad
 \norm u_\pi=
 \inf_{u=\sum_i a_i\otimes b_i}\sum_i\norm{a_i}\norm{b_i}.
\]
The completions are denoted by $E\te F$ and $E\tp F$. The tensor $u$
defines an operator
\[
 E^*\longrightarrow F,\qquad x^*\longmapsto
 \sum_{i=1}^q x^*(x_i)y_i,
\]
whose operator norm equals $\norm u_\eps$. In finite dimensions this
identifies $E\otimes_\eps F$ with all operators from $E^*$ to $F$.
Moreover, $(E\tp F)^*$ is the space of bounded bilinear forms on
$E\times F$, with its usual norm; see \cite{Ryan}.

Let $H=(\R^2,\norm\cdot_2)$, with standard basis $e_1,e_2$, and identify
$H^*$ with $H$ by the inner product. In particular, the map
\begin{equation}\label{eq:matrix}
 H\otimes_\eps H\longrightarrow (M_2(\R),\norm\cdot_{\mathrm{op}}),
 \qquad x\otimes y\longmapsto yx^{\mathsf T},
\end{equation}
is an onto linear isometry. Under this identification
\[
 z_0=e_1\otimes e_1+e_2\otimes e_2
\]
corresponds to $\id_H$ and has injective norm one. We put $\ell_1^0=\ell_\infty^0=\{0\}$, and interpret empty
sums as zero.

\begin{theorem}\label{thm:main}
There exist separable Banach spaces $X$ and $Y$, both with the Daugavet
property, such that neither
\[
 Z_\eps=X\te X\qquad\text{nor}\qquad Z_\pi=Y\tp Y
\]
has the Daugavet property. More precisely, $B_{Z_\eps}$ has a strongly
exposed point, and $B_{Z_\pi^*}$ has a weak-star strongly exposed point.
There are rank-one operators $T_\eps$ on $Z_\eps$ and $T_\pi$ on $Z_\pi$
such that, for $\nu\in\{\eps,\pi\}$,
\[
 \norm{T_\nu}=1,\qquad
 \norm{\id_{Z_\nu}+T_\nu}\le\frac{199}{100}<2.
\]
The space $X$ contains $H$ isometrically as the range of a contractive
projection, and $Y$ is a closed subspace of $X^*$.
\end{theorem}

In Section~2 we describe the finite extension of an arbitrary finite-dimensional space $E$. In Section~3 we extend $H$ infinitely many times to construct
$X$. In Section~4, we show that the norms on the finite-dimensional spaces arising in our construction can be recovered from a family of linear maps into spaces of a particularly simple form. Section~5 then establishes a quantitative estimate for operators on these spaces, with constants independent of their dimensions. The injective example
is completed in Section~6. In Sections~7 and~8 we construct the projective example $Y$. Section~9 records the relation with the operator Daugavet
properties. Section 10 is independent from the rest and presents a result on the isomorphic theory of Daugavet spaces. 

\section{A finite-dimensional extension}

Let $E$ be a finite-dimensional normed space, let $p\in S_E$, and let
$m\ge2$ be an integer. We define
\begin{equation}\label{eq:G}
 G(E,p,m)=E\oplus\R^m,\qquad
 \norm{(x,a)}_G=
 \max_{1\le j\le m}\bigl(\norm{x-a_jp}_E+|a_j|\bigr).
\end{equation}
We always identify $E$ with $E\oplus\{0\}$.

\begin{lemma}\label{lem:extension}
Formula \eqref{eq:G} defines a norm. The inclusion $E\hookrightarrow G$
is isometric, and the maps
\[
 R:G\to E,\quad R(x,a)=x,\qquad
 A_j:G\to E\oplus_1\R,\quad A_j(x,a)=(x-a_jp,a_j)
\]
are contractions. If $u_1,\ldots,u_m$ is the standard basis of $\R^m$
and $h_j=(p,u_j)$, then
\begin{equation}\label{eq:extensionproperties}
 \norm{h_j}=1,\qquad
 \norm{x+h_j}=2\quad(x\in S_E),\qquad
 \norm{\frac1m\sum_{j=1}^m h_j-p}=\frac2m.
\end{equation}
\end{lemma}

\begin{proof}
Homogeneity and the triangle inequality follow from the corresponding
properties of the expressions in the maximum. If the maximum is zero,
then every $a_j$ is zero, and consequently $x=0$. For $a=0$ the maximum
equals $\norm x_E$. Each $A_j$ is contractive by the definition, and
\[
 \norm x_E\le\norm{x-a_jp}_E+|a_j|\le\norm{(x,a)}_G
\]
proves the assertion about $R$.

For $h_j$, all expressions in the maximum equal one. For $x+h_j$, the
$j$th expression equals $\norm x_E+1=2$, and the triangle inequality
gives the reverse bound. Finally, the difference in the last part of
\eqref{eq:extensionproperties} is $(0,m^{-1},\ldots,m^{-1})$, whose norm
is $2/m$.
\end{proof}

\section{The inductive construction}

Put $E_0=H$ and $N_{-1}=\varnothing$. Suppose that $E_n$ has been
constructed. Choose a finite set
\[
 N_n=\{p_{n,1},\ldots,p_{n,r_n}\}\subset S_{E_n},\qquad
 N_{n-1}\subset N_n,
\]
such that every point of $S_{E_n}$ has distance strictly less than
$\eps_n=2^{-n}$ from $N_n$. This is possible: first take a finite open
$\eps_n$-net, and then add all points of $N_{n-1}$. These points still
have norm one because the earlier inclusions are isometric. Set
$m_n=2^{n+3}$ and define successively
\begin{equation}\label{eq:outerstages}
 E_{n,0}=E_n,\qquad
 E_{n,j}=G(E_{n,j-1},p_{n,j},m_n)\quad(1\le j\le r_n),\qquad
 E_{n+1}=E_{n,r_n}.
\end{equation}
In particular, all $p_{n,j}$ are chosen in the same space $E_n$, before
the extensions in \eqref{eq:outerstages} are performed. Put
\[
 X=\overline{\bigcup_{n\ge0}E_n}.
\]

It will be convenient to number the elementary extensions as well.
Let $b_0=0$ and $b_{n+1}=b_n+r_n$. For $0\le j\le r_n$ set
$V_{b_n+j}=E_{n,j}$. Thus $V_0=H$, $E_n=V_{b_n}$, and
\[
 V_{s+1}=G(V_s,p_s,\mu_s),\qquad
 p_{b_n+j-1}=p_{n,j},\quad \mu_{b_n+j-1}=m_n.
\]
Here the notation $p_s$ records an elementary extension point; it does
not denote a new choice of a net.

\begin{lemma}\label{lem:projections}
For every $s\ge0$ there is a contractive projection $\rho_s:X\to V_s$
such that, when $s\le t$,
\[
 \rho_s\rho_t=\rho_s,\qquad \rho_t\rho_s=\rho_s.
\]
Moreover, $\rho_sx\to x$ for every $x\in X$.
\end{lemma}

\begin{proof}
For each elementary extension let
\[
 R_s:V_{s+1}\longrightarrow V_s,\qquad R_s(x,a)=x.
\]
By Lemma~\ref{lem:extension}, $R_s$ is contractive and fixes $V_s$.
For integers $t>s$ define the map from $V_t$ to $V_s$ by
\[
 R_{s,t}=R_s\circ R_{s+1}\circ\cdots\circ R_{t-1},
 \qquad R_{s,s}=\id_{V_s}.
\]
The rightmost map $R_{t-1}$ is applied first. Thus $R_{s,t}$ is a
contraction, and for $s\le t\le u$ one has
\[
 R_{s,u}=R_{s,t}\circ R_{t,u},\qquad
 R_{s,u}|_{V_t}=R_{s,t}.
\]
The second equality follows because $R_{t,u}$ fixes every vector of
$V_t$. Consequently the prescriptions
$\rho_sx=R_{s,t}x$ for $x\in V_t$, $t\ge s$, agree on intersections.
They define a contraction on $\bigcup_tV_t$ and hence extend uniquely to
a contraction $\rho_s:X\to V_s$. Its restriction to $V_s$ is the
identity, so it is a projection. The displayed compatibility relations
first hold on each $V_u$ and then on $X$ by continuity.

Given $x\in X$ and $\eta>0$, choose $v\in V_t$ with
$\norm{x-v}<\eta/2$. For $s\ge t$ we have $\rho_sv=v$, and hence
\[
 \norm{\rho_sx-x}
 \le\norm{\rho_s(x-v)}+\norm{v-x}
 \le2\norm{x-v}<\eta.
\]
This proves the convergence.
\end{proof}

Write $P=\rho_0:X\to H$ and $P_s=P|_{V_s}$. We will repeatedly use
$P\rho_s=P$. The nested choice of the nets gives the following two
facts:
\begin{equation}\label{eq:recurrence}
 \mathcal P=\bigcup_{n\ge0}N_n\text{ is dense in }S_X,
 \quad\text{and every }p\in\mathcal P\text{ occurs as }p_s
 \text{ for arbitrarily large }s.
\end{equation}
For the first assertion, choose a normalized vector $v\in E_N$ close
to a given point of $S_X$. For any $n\ge N$, this same vector belongs
to $S_{E_n}$, and some $p\in N_n$ satisfies $\norm{v-p}<\eps_n$.
Letting $n$ increase proves the assertion. The second follows from
$N_n\subset N_{n+1}$ and the fact that each point of $N_n$ is used in
one of the extensions leading from $E_n$ to $E_{n+1}$.

\begin{proposition}\label{prop:DPX}
The separable Banach space $X$ has the Daugavet property.
\end{proposition}

\begin{proof}
Fix $f,g\in S_X$ and $\eta>0$. Choose $0<\tau<\eta/3$, and then choose
$N$ and $f_0,g_0\in S_{E_N}$ such that
$\norm{f-f_0}<\tau$ and $\norm{g-g_0}<\tau$. Normalized approximants
exist because $\bigcup_nE_n$ is dense in $X$. Take $n\ge N$ so large
that
\[
 \eps_n+\frac{2}{m_n}<\tau.
\]
The vector $f_0$ belongs to $S_{E_n}$. By the defining net property,
there is $j\in\{1,\ldots,r_n\}$ with
$\norm{f_0-p_{n,j}}<\eps_n$. Apply
Lemma~\ref{lem:extension} at the extension
$E_{n,j}=G(E_{n,j-1},p_{n,j},m_n)$, and denote its new vectors by
$h_1,\ldots,h_{m_n}$. Both $g_0$ and $p_{n,j}$ belong to the preceding
space $E_{n,j-1}$, and retain norm one there. Therefore
\[
 \norm{h_k}=1,\quad \norm{g_0+h_k}=2,
 \quad
 \norm{\frac1{m_n}\sum_{k=1}^{m_n}h_k-p_{n,j}}=\frac2{m_n}.
\]
These equalities remain true in $X$. It follows that
\begin{align*}
 \norm{\frac1{m_n}\sum_{k=1}^{m_n}h_k-f}
 &\le \norm{\frac1{m_n}\sum_{k=1}^{m_n}h_k-p_{n,j}}
      +\norm{p_{n,j}-f_0}+\norm{f_0-f}\\
 &<\frac2{m_n}+\eps_n+\tau<2\tau<\eta,\\
 \norm{g+h_k}
 &\ge\norm{g_0+h_k}-\norm{g-g_0}>2-\tau>2-\eta.
\end{align*}

The Daugavet property follows now from Theorem 1.1 (3).
\end{proof}

\section{Norming representations of the finite stages}

For $r\ge0$ put $F_r=H\oplus_1\ell_1^r$. If
$K:\ell_1^r\to H$ is a contraction, then so is
\[
 R_K:F_r\to H,\qquad R_K(v,b)=v+Kb,
\]
since $\norm{v+Kb}_2\le\norm v_2+\norm b_1$.

\begin{lemma}\label{lem:norming}
For each $s\ge0$ there are a finite nonempty index set $\Gamma_s$,
contractions $J_\gamma:V_s\to F_s$, $\gamma\in\Gamma_s$, and a
contraction $K_s:\ell_1^s\to H$, independent of $\gamma$, such that
\begin{equation}\label{eq:norming}
 \norm x=\max_{\gamma\in\Gamma_s}\norm{J_\gamma x},\qquad
 J_\gamma h=(h,0)\ (h\in H),\qquad
 R_{K_s}J_\gamma=P_s.
\end{equation}
\end{lemma}

\begin{proof}
At $s=0$ take a singleton index set, let $K_0:\{0\}\to H$ be zero,
and define
\[
 J_\gamma:H\longrightarrow H\oplus_1\{0\},\qquad
 J_\gamma h=(h,0).
\]
All assertions in \eqref{eq:norming} then hold.

Suppose they hold for $V_s$, and write
$V_{s+1}=G(V_s,p_s,\mu_s)$. Set
$\Gamma_{s+1}=\Gamma_s\times\{1,\ldots,\mu_s\}$. Under the natural
identification $F_{s+1}=F_s\oplus_1\R$, define
\[
 J'_{\gamma,j}(x,a)=\bigl(J_\gamma(x-a_jp_s),a_j\bigr),
 \qquad
 K_{s+1}(b,t)=K_sb+tP_sp_s.
\]
Since $\norm{P_sp_s}\le1$, one has
$\norm{K_{s+1}(b,t)}\le\norm b_1+|t|$. Moreover,
\begin{align*}
 \max_{\gamma,j}\norm{J'_{\gamma,j}(x,a)}
 &=\max_j\left(\max_\gamma\norm{J_\gamma(x-a_jp_s)}+|a_j|\right)\\
 &=\max_j\bigl(\norm{x-a_jp_s}+|a_j|\bigr)
 =\norm{(x,a)}_{V_{s+1}}.
\end{align*}
In particular each $J'_{\gamma,j}$ is contractive. For $h\in H$, its
new coordinates vanish, so $J'_{\gamma,j}h=((h,0),0)=(h,0)$.
Finally, if $J_\gamma(x-a_jp_s)=(v,b)$, then
\begin{align*}
 R_{K_{s+1}}J'_{\gamma,j}(x,a)
 &=v+K_sb+a_jP_sp_s\\
 &=P_s(x-a_jp_s)+a_jP_sp_s
 =P_sx=P_{s+1}(x,a).
\end{align*}
This proves the induction. Notice that $K_{s+1}$ does not depend on
either $\gamma$ or $j$.
\end{proof}

\begin{lemma}\label{lem:tensornorming}
Suppose that $J_\gamma:E\to U_\gamma$ and $L_\eta:F\to W_\eta$ are
finite families of linear maps satisfying
\[
 \norm x=\max_\gamma\norm{J_\gamma x},\qquad
 \norm y=\max_\eta\norm{L_\eta y}.
\]
Then, for every $u\in E\otimes F$,
\begin{equation}\label{eq:tensornorming}
 \norm u_\eps=\max_{\gamma,\eta}
       \norm{(J_\gamma\otimes L_\eta)u}_\eps.
\end{equation}
\end{lemma}

\begin{proof}
Let $\mathcal U=(\bigoplus_\gamma U_\gamma)_\infty$ and
$\mathcal W=(\bigoplus_\eta W_\eta)_\infty$. The maps
$Jx=(J_\gamma x)_\gamma$ and $Ly=(L_\eta y)_\eta$ are linear
isometries.

Let $u=\sum_{i=1}^n x_i\otimes y_i\in E\otimes F$.
By the definition of the injective tensor norm,
\[
\begin{aligned}
\lVert (J\otimes L)u\rVert_\varepsilon
&=\sup_{\substack{\varphi\in B_{\mathcal U^*}\\
                   \psi\in B_{\mathcal W^*}}}
\left|\sum_{i=1}^n
\varphi(Jx_i)\psi(Ly_i)\right|.
\end{aligned}
\]
Since $J$ and $L$ are isometries, the Hahn--Banach
theorem gives
\[
J^*(B_{\mathcal U^*})=B_{E^*},
\qquad
L^*(B_{\mathcal W^*})=B_{F^*}.
\]
Indeed, every $f\in B_{E^*}$ induces a functional
on $J(E)$ of the same norm by $Jx\mapsto f(x)$,
and this functional extends to $\mathcal U$
without increasing its norm.
The same argument applies to $L$.

Consequently,
\[
\begin{aligned}
\lVert (J\otimes L)u\rVert_\varepsilon
&=\sup_{\substack{f\in B_{E^*}\\
                  g\in B_{F^*}}}
\left|\sum_{i=1}^n f(x_i)g(y_i)\right|\\
&=\lVert u\rVert_\varepsilon.
\end{aligned}
\]

It remains to compute the tensor norm in the two finite sums. Their
duals are the corresponding $\ell_1$-sums. If
$v =\sum_{i=1}^n x_i\otimes y_i\in\mathcal U\otimes\mathcal W$ has components $v_{\gamma,\eta} = \sum_{i=1}^n x_{i,\gamma}\otimes y_{i,\eta}$ for $x_i = (x_{i,\gamma})_\gamma$ and $y_i = (y_{i,\eta})_\eta$,
then functionals $a^*=(a^*_\gamma)$ and $b^*=(b^*_\eta)$ with
$\sum_\gamma\norm{a^*_\gamma}\le1$ and
$\sum_\eta\norm{b^*_\eta}\le1$ satisfy
\[
\begin{aligned}
(a^*\otimes b^*)(v)
&=\sum_{i=1}^n a^*(x_i)b^*(y_i)\\
&=\sum_{i=1}^n
  \left(\sum_\gamma
  a_\gamma^*(x_{i,\gamma})\right)
  \left(\sum_\eta
  b_\eta^*(y_{i,\eta})\right)\\
&=\sum_{\gamma,\eta}
  (a_\gamma^*\otimes b_\eta^*)
  (v_{\gamma,\eta}).
\end{aligned}
\]

By the triangle inequality and the definition
of the injective tensor norm,
\[
\begin{aligned}
|(a^*\otimes b^*)(v)|
&\le
\sum_{\gamma,\eta}
\left|
(a_\gamma^*\otimes b_\eta^*)
(v_{\gamma,\eta})
\right|\\
&\le
\sum_{\gamma,\eta}
\|a_\gamma^*\|
\|b_\eta^*\|
\|v_{\gamma,\eta}\|_\varepsilon\\
&\le
\left(
\max_{\gamma,\eta}
\|v_{\gamma,\eta}\|_\varepsilon
\right)
\left(\sum_\gamma\|a_\gamma^*\|\right)
\left(\sum_\eta\|b_\eta^*\|\right)\\
&\le
\max_{\gamma,\eta}
\|v_{\gamma,\eta}\|_\varepsilon.
\end{aligned}
\]

Taking the supremum over all
$a^*\in B_{\mathcal U^*}$ and
$b^*\in B_{\mathcal W^*}$ gives
\[
\|v\|_\varepsilon
\le
\max_{\gamma,\eta}
\|v_{\gamma,\eta}\|_\varepsilon.
\]

For the reverse inequality, fix one pair $(\gamma,\eta)$ and choose
functionals supported only in these coordinates, with evaluations
arbitrarily close to $\norm{v_{\gamma,\eta}}_\eps$. Thus
$\norm v_\eps=\max_{\gamma,\eta}\norm{v_{\gamma,\eta}}_\eps$.
Apply this identity to $v=(J\otimes L)u$.
\end{proof}

\section{A uniform operator estimate}

\begin{lemma}\label{lem:block}
Let $U=H\oplus_1\ell_1^k$ and $V=H\oplus_1\ell_1^l$, and let
$K:\ell_1^k\to H$ and $L:\ell_1^l\to H$ be contractions. Suppose
that $W:U^*\to V$ is linear and $\norm W\le1$. Put
\[
 M=R_LWR_K^*:H\to H,\qquad
 \delta=1-\tfrac12\tr M,\qquad
 I_0(h,t)=(h,0)\quad((h,t)\in U^*).
\]
Then $0\le\delta\le2$ and
\begin{equation}\label{eq:blockbound}
 \norm{W-I_0}\le2\sqrt\delta+122\delta.
\end{equation}
\end{lemma}

\begin{proof}
Using $U^*=H\oplus_\infty\ell_\infty^k$, write
\[
 W(h,t)=(Ah+Bt,Ch+Dt),
\]
where the four bounded linear operators are
\[
A:H\to H,\quad
B:\ell_\infty^k\to H,\quad
C:H\to\ell_1^l,\quad
D:\ell_\infty^k\to\ell_1^l.
\]
The contraction assumption means that
\begin{equation}\label{eq:contraction}
 \norm{Ah+Bt}_2+\norm{Ch+Dt}_1\le1
 \quad(\norm h_2\le1,\ \norm t_\infty\le1).
\end{equation}
Since $R_K$ and $R_L$ are contractions, $\norm M\le1$. In particular,
$-2\le\tr M\le2$, proving the assertion about $\delta$.

Set
\begin{equation}\label{eq:AN}
 A_0=A+BK^*,\qquad N=C+DK^*.
\end{equation}
Then $R_K^*h=(h,K^*h)$ and
\begin{equation}\label{eq:Mdecomp}
 M=A_0+LN,\qquad
 \norm{A_0h}_2+\norm{Nh}_1\le1\quad(h\in S_H).
\end{equation}

Indeed, for $h \in H$:

\begin{align*}
Mh
&=R_LWR_K^*h
=R_LW(h,K^*h)\\
&=R_L(Ah+BK^*h,\;Ch+DK^*h)\\
&=Ah+BK^*h+L(Ch+DK^*h)\\
&=(A+BK^*)h+L(C+DK^*)h\\
&=A_0h+LNh.
\end{align*}

Moreover, if $h\in S_H$, then
\[
\|K^*h\|_\infty
\le \|K^*\|\|h\|_2
\le 1.
\]
Thus $(h,K^*h)$ belongs to the unit ball of
$U^*=H\oplus_\infty\ell_\infty^k$.
Applying the contraction estimate (5.2) with
$t=K^*h$, we obtain
\[
\begin{aligned}
1
&\ge \|Ah+BK^*h\|_2
   +\|Ch+DK^*h\|_1\\
&=\|A_0h\|_2+\|Nh\|_1.
\end{aligned}
\]

Let $\av$ denote integration with respect to normalized arc length on
$S_H$. Parametrizing $h=(\cos\theta,\sin\theta)$ gives
\begin{equation}\label{eq:circle}
 \av(hh^{\mathsf T})=\tfrac12\id_H,
 \qquad
 \av|\ip v h|=\frac2\pi\norm v_2\quad(v\in H).
\end{equation}

To verify these identities, write
\[
h(t)=\cos t\,e_1+\sin t\,e_2,
\qquad 0\le t<2\pi.
\]
Recall that integration with respect to the
normalized uniform measure on $S_H$ is given by
\[
\av f(h)
=\frac{1}{2\pi}\int_0^{2\pi}f(h(t))\,dt.
\]

For the first identity, we have
\[
\begin{aligned}
h(t)h(t)^{\mathsf T}
&=\cos^2t\,e_1e_1^{\mathsf T}
 +\sin^2t\,e_2e_2^{\mathsf T}+\cos t\sin t\,
 \bigl(e_1e_2^{\mathsf T}
      +e_2e_1^{\mathsf T}\bigr).
\end{aligned}
\]
Since
\begin{align*}
\int_0^{2\pi}\cos^2t\,dt
&=\int_0^{2\pi}\sin^2t\,dt=\pi,
\end{align*}
and $\int_0^{2\pi}\cos t\sin t\,dt=0$, it follows that
\[
\begin{aligned}
\av(hh^{\mathsf T})
&=\frac{1}{2\pi}
  \int_0^{2\pi}h(t)h(t)^{\mathsf T}\,dt\\
&=\frac12
  \bigl(e_1e_1^{\mathsf T}
       +e_2e_2^{\mathsf T}\bigr)\\
&=\frac12\id_H.
\end{aligned}
\]

For the second identity, the case $v=0$ is
immediate. If $v\ne0$, choose $\theta\in\mathbb R$
such that
\[
v=\norm v_2
  (\cos\theta\,e_1+\sin\theta\,e_2).
\]
Then
\[
\begin{aligned}
\av|\ip v h|
&=\frac{1}{2\pi}
  \int_0^{2\pi}|\ip v{h(t)}|\,dt\\
&=\frac{\norm v_2}{2\pi}
  \int_0^{2\pi}
  |\cos\theta\cos t+\sin\theta\sin t|\,dt\\
&=\frac{\norm v_2}{2\pi}
  \int_0^{2\pi}|\cos(t-\theta)|\,dt\\
&=\frac{\norm v_2}{2\pi}
  \int_0^{2\pi}|\cos t|\,dt\\
&=\frac{4}{2\pi}\norm v_2\\
&=\frac2\pi\norm v_2.
\end{aligned}
\]
Here we used the periodicity of $|\cos t|$
and the identity
\[
\int_0^{2\pi}|\cos t|\,dt
=4\int_0^{\pi/2}\cos t\,dt=4.
\]

Write $(Nh)_r=\ip{n_r}h$ for $n_r\in H$, put $l_r=Lu_r$, and set
$S=\sum_{r=1}^l\norm{n_r}_2$. The sum is finite, so
\[
 \av\norm{Nh}_1
 =\sum_{r=1}^l\av|\ip{n_r}h|
 =\frac2\pi\sum_{r=1}^l\norm{n_r}_2
 =\frac2\pi S.
\]
Integrating \eqref{eq:Mdecomp} therefore yields
\begin{equation}\label{eq:averageupper}
 \av\norm{A_0h}_2+\frac2\pi S\le1.
\end{equation}
On the other hand, the left equation in \eqref{eq:circle} implies
$\av\ip h{Mh}=\tfrac12\tr M$. Indeed,
\[
\langle h,Mh\rangle
=\operatorname{tr}(Mhh^{\mathsf T})
\]
is elementary, and integration commutes with
both $M$ and the trace. Since
$LNh=\sum_{r=1}^l\ip{n_r}h\,l_r$, we obtain together with (5.4)
\begin{align}
 1-\delta
 &= \tfrac12\tr M = \av\ip h{Mh} =\av\ip h{A_0h}
   +\sum_{r=1}^l\av\bigl(\ip{n_r}h\ip h{l_r}\bigr)\notag\\
 &=\av\ip h{A_0h}+\frac12\sum_{r=1}^l\ip{n_r}{l_r}\notag\le\av\norm{A_0h}_2+\frac12\sum_{r=1}^l\norm{n_r}_2\norm{l_r}_2 \\
 &\le\av\norm{A_0h}_2+\frac12S.\label{eq:averagelower}
\end{align}
Here $\norm h_2=1$ and $\norm{l_r}_2=\norm{Lu_r}_2\le1$.
Combining \eqref{eq:averageupper} and \eqref{eq:averagelower}, we get
\[
 1-\delta\le1-\left(\frac2\pi-\frac12\right)S.
\]
As $2/\pi-1/2>1/8$, it follows that $S\le8\delta$.
Also \eqref{eq:Mdecomp} gives $\norm{A_0h}_2\le1$,
whereas \eqref{eq:averagelower} gives
$\av\norm{A_0h}_2\ge1-\delta-S/2$. Consequently
\begin{equation}\label{eq:deficit}
 S\le8\delta,\qquad
 0\le1-\av\norm{A_0h}_2
 \le\delta+\frac12S\le5\delta.
\end{equation}
In addition, $\norm N_{2\to1}\le S$, because
$\sum_r|\ip{n_r}h|\le S\norm h_2$.

We next estimate $D$. Assume first that $k,l>0$. Denote the columns
of $K$ by $k_i=Ku_i$, and put
\[
 a=1-\frac2\pi>\frac13,\qquad
 d_i=\av\bigl(1-|\ip{k_i}h|\bigr)
     =1-\frac2\pi\norm{k_i}_2\ge a.
\]
For a fixed $t\in B_{\ell_\infty^k}$ and $h\in S_H$, define
\[
 (t_h)_i=\frac a{d_i}t_i\bigl(1-|\ip{k_i}h|\bigr).
\]
For each coordinate,
\[
 |(K^*h)_i\pm(t_h)_i|
 \le|\ip{k_i}h|+\frac a{d_i}|t_i|
            \bigl(1-|\ip{k_i}h|\bigr)\le1.
\]
Thus $K^*h\pm t_h\in B_{\ell_\infty^k}$. The definition of $d_i$
also gives $\av t_h=at$ coordinate by coordinate.
Applying \eqref{eq:contraction} with $K^*h+t_h$ and $K^*h-t_h$ gives
\[
 \norm{A_0h+Bt_h}_2+\norm{Nh+Dt_h}_1\le1,
 \qquad
 \norm{A_0h-Bt_h}_2+\norm{Nh-Dt_h}_1\le1.
\]
The sum of the two Hilbert-space terms is at least
$2\norm{A_0h}_2$, and
$\norm{Nh+Dt_h}_1+\norm{Nh-Dt_h}_1\ge2\norm{Dt_h}_1$.
Adding the inequalities and dividing by two therefore yields
\[
 \norm{A_0h}_2+\norm{Dt_h}_1\le1.
\]
By linearity, the triangle inequality for the integral, and
\eqref{eq:deficit},
\[
 a\norm{Dt}_1
 =\norm{D(\av t_h)}_1
 =\norm{\av(Dt_h)}_1
 \le\av\norm{Dt_h}_1
 \le1-\av\norm{A_0h}_2\le5\delta.
\]
Taking the supremum over $t\in B_{\ell_\infty^k}$ gives
$\norm D_{\infty\to1}\le15\delta$. If $k=0$ or $l=0$, the same
bound holds because $D=0$. Since $C=N-DK^*$ and $\norm{K^*}\le1$,
we have proved
\begin{equation}\label{eq:CD}
 \norm D_{\infty\to1}\le15\delta,
 \qquad
 \norm C_{2\to1}\le\norm N_{2\to1}+\norm D_{\infty\to1}
 \le23\delta.
\end{equation}

All spaces in this lemma are finite-dimensional. The adjoint operator
$W^*:V^*\to U$ is consequently a contraction of the same form, with
block matrix
\[
 W^*=\begin{pmatrix}A^*&C^*\\ B^*&D^*\end{pmatrix}.
\]
Interchanging $K$ and $L$, its compressed operator is
$R_KW^*R_L^*=M^*$, whose trace equals $\tr M$. Apply the bound for
the lower left block in \eqref{eq:CD} to this operator. Since that
block is $B^*:H\to\ell_1^k$, we obtain
\begin{equation}\label{eq:B}
 \norm B_{\infty\to2}=\norm{B^*}_{2\to1}\le23\delta.
\end{equation}

Let $\norm\cdot_{\mathrm{HS}}$ be the Hilbert--Schmidt norm on
operators on $H$. Since $\norm M\le1$ and $\dim H=2$,

\begin{align*}
\|M-\id_H\|^2
&\le \|M-\id_H\|_{\mathrm{HS}}^2
=\sum_{i=1}^2\|(M-\id_H)e_i\|_2^2\\
&=\sum_{i=1}^2\|Me_i-e_i\|_2^2
=\sum_{i=1}^2
\left(\|Me_i\|_2^2+\|e_i\|_2^2
-2\langle Me_i,e_i\rangle\right)\\
&\le \sum_{i=1}^2
\left(2-2\langle Me_i,e_i\rangle\right)=4-2\sum_{i=1}^2\langle Me_i,e_i\rangle\\
&=4-2\tr M=4-4(1-\delta)=4\delta.
\end{align*}
Equations \eqref{eq:AN} and \eqref{eq:Mdecomp} give
$A=M-BK^*-LC-LDK^*$. Hence \eqref{eq:CD} and \eqref{eq:B} yield
\[
 \norm{A-\id_H}\le2\sqrt\delta+(23+23+15)\delta
 =2\sqrt\delta+61\delta.
\]
Finally, for $\norm h_2,\norm t_\infty\le1$,
\begin{align*}
 \norm{(W-I_0)(h,t)}
 &=\norm{(A-\id_H)h+Bt}_2+\norm{Ch+Dt}_1\\
 &\le\norm{A-\id_H}+\norm B+\norm C+\norm D\\
 &\le2\sqrt\delta+122\delta.
\end{align*}
This proves \eqref{eq:blockbound}.
\end{proof}

\section{The injective tensor square}

Recall that $P_s=P|_{V_s}$ and $z_0=e_1\otimes e_1+e_2\otimes e_2$.
For $w\in V_s\otimes V_s$ put
\[
 \Phi_s(w)=\tfrac12\tr\bigl((P_s\otimes P_s)w\bigr),
\]
using the matrix identification \eqref{eq:matrix}.

\begin{proposition}\label{prop:finitebound}
For every $s\ge0$ and every $w\in V_s\otimes V_s$ with
$\norm w_\eps\le1$,
\begin{equation}\label{eq:finitebound}
 \norm{w-z_0}_\eps
 \le2\sqrt{1-\Phi_s(w)}+122\bigl(1-\Phi_s(w)\bigr).
\end{equation}
\end{proposition}

\begin{proof}
Apply Lemma~\ref{lem:norming} to $V_s$. For a fixed pair
$\gamma,\eta\in\Gamma_s$, regard
\[
 W_{\gamma,\eta}=(J_\gamma\otimes J_\eta)w
\]
as an operator $F_s^*\to F_s$. Its norm is at most one because
$J_\gamma$ and $J_\eta$ are contractions. We apply
Lemma~\ref{lem:block} with
\[
 U=V=F_s,\qquad k=l=s,\qquad K=L=K_s,
 \qquad W=W_{\gamma,\eta}.
\]
By the factorization $R_{K_s}J_\gamma=P_s$ in
Lemma~\ref{lem:norming}, the corresponding compressed operator is
\begin{align*}
 R_{K_s}W_{\gamma,\eta}R_{K_s}^*
 &=(R_{K_s}J_\gamma\otimes R_{K_s}J_\eta)w\\
 &=(P_s\otimes P_s)w.
\end{align*}
The first equality is read using the tensor-to-operator identification:
applying maps to the first and second tensor factors corresponds to
composition on the right by the adjoint and on the left, respectively.
Thus every pair has the same trace defect
$\delta=1-\Phi_s(w)$.

The identities $J_\gamma e_i=(e_i,0)$ imply that
$(J_\gamma\otimes J_\eta)z_0$ is the operator
$I_0(h,t)=(h,0)$. Lemma~\ref{lem:block} therefore gives, for every pair,
\[
 \norm{(J_\gamma\otimes J_\eta)(w-z_0)}_\eps
 \le2\sqrt\delta+122\delta.
\]
Taking the maximum over $\gamma,\eta$ and using
Lemma~\ref{lem:tensornorming} proves \eqref{eq:finitebound}.
\end{proof}

Put $Z_\eps=X\te X$ and define
\begin{equation}\label{eq:Phi}
 \Phi(z)=\tfrac12\tr\bigl((P\otimes P)z\bigr)
 \qquad(z\in Z_\eps).
\end{equation}
The map $P\otimes P$ is contractive, and
$|\tr A|/2\le\norm A_{\mathrm{op}}$ on $M_2(\R)$. Moreover,
$\norm{z_0}_\eps=1$ by the injectivity of the tensor norm and
\eqref{eq:matrix}, and $\Phi(z_0)=1$. Consequently $\norm\Phi=1$.

\begin{proposition}\label{prop:globalbound}
For every $z\in B_{Z_\eps}$,
\begin{equation}\label{eq:globalbound}
 \norm{z-z_0}_\eps
 \le2\sqrt{1-\Phi(z)}+122\bigl(1-\Phi(z)\bigr).
\end{equation}
In particular, $\Phi$ strongly exposes $z_0$ in $B_{Z_\eps}$.
\end{proposition}

\begin{proof}
By Lemma~\ref{lem:projections}, the operators
$Q_s=\rho_s\otimes\rho_s$ extend to contractions on $Z_\eps$, with
range in $V_s\otimes V_s$. For an elementary tensor,
\[
 \norm{Q_s(x\otimes y)-x\otimes y}_\eps
 \le\norm{\rho_sx-x}\norm{\rho_sy}
     +\norm x\norm{\rho_sy-y}\longrightarrow0.
\]
The same holds for finite sums. By density and $\norm{Q_s}\le1$,
it holds for every $z\in Z_\eps$. Furthermore, $P\rho_s=P$, so
\[
 \Phi_s(Q_sz)=\Phi(z).
\]
For $z\in B_{Z_\eps}$, Proposition~\ref{prop:finitebound} applies to
$Q_sz$ without any normalization and gives
\[
 \norm{Q_sz-z_0}_\eps
 \le2\sqrt{1-\Phi(z)}+122\bigl(1-\Phi(z)\bigr).
\]
Passing to the limit proves \eqref{eq:globalbound}. Its right-hand
side tends to zero when $\Phi(z)$ tends to one, which is precisely
the asserted strong exposure.
\end{proof}

\begin{corollary}\label{cor:injective}
The operator $T_\eps z=-\Phi(z)z_0$ on $Z_\eps$ has rank one,
$\norm{T_\eps}=1$, and
$\norm{\id_{Z_\eps}+T_\eps}\le199/100$.
\end{corollary}

\begin{proof}
The norm and rank assertions follow from
$\norm{z_0}=\norm\Phi=\Phi(z_0)=1$. Let $z\in B_{Z_\eps}$.
If $|\Phi(z)|\le99/100$, then
\[
 \norm{z-\Phi(z)z_0}\le1+|\Phi(z)|\le199/100.
\]
Otherwise let $w=\operatorname{sign}(\Phi(z))z$, so that
$\Phi(w)=|\Phi(z)|>99/100$. Proposition~\ref{prop:globalbound} gives
\[
 \norm{w-z_0}\le2\sqrt{1/100}+122/100=71/50.
\]
Therefore
\[
 \norm{z-\Phi(z)z_0}
 =\norm{w-\Phi(w)z_0}
 \le\norm{w-z_0}+1-\Phi(w)\le143/100.
\]
Taking the supremum over $z$ proves the result.
\end{proof}

\section{The coordinate dual space}

For each $s\ge0$, define
\begin{equation}\label{eq:Y}
 Y_s=\rho_s^*(V_s^*)\subset X^*,\qquad
 Y=\overline{\bigcup_{s\ge0}Y_s}^{\norm\cdot_{X^*}}.
\end{equation}
The adjoint map $\rho_s^*:V_s^*\to X^*$ is an isometry. Indeed,
$\norm{\rho_s^*f}\le\norm f$ by contractivity, while evaluating on
$B_{V_s}\subset B_X$ gives the reverse inequality. The compatibility
of the projections also gives $Y_s\subset Y_t$ for $s\le t$:
an element $f\circ\rho_s$ equals
$(f\circ R_{s,t})\circ\rho_t$. In particular, $Y$ is separable.
An element of $Y_s$ vanishes on the new coordinate directions of every
later extension. We use the inherited norm on $Y$; no identification
of $Y^*$ with $X$ is required.

\begin{lemma}\label{lem:coordinates}
Let $p\in S_{V_M}$ and suppose that
$M\le s_1<\cdots<s_r$ satisfy $p_{s_k}=p$. At the extension
$V_{s_k+1}=G(V_{s_k},p,\mu_{s_k})$, let
$c_k:V_{s_k+1}\to\R$, $c_k(x,a)=a_1$, and put
$\lambda_k=c_k\circ\rho_{s_k+1}\in Y_{s_k+1}$. Then
\begin{equation}\label{eq:coordinates}
 \sum_{k=1}^r|\lambda_k(x)|\le\norm x\quad(x\in X).
\end{equation}
In particular, $\norm{\sum_{k=1}^r\lambda_k}\le1$.
\end{lemma}

\begin{proof}
We construct contractions
\[
 J_t:V_t\longrightarrow V_M\oplus_1\ell_1^{d(t)},
 \qquad d(t)=\#\{k:s_k<t\},\qquad M\le t\le s_r+1,
\]
such that $J_tv=(v,0)$ for $v\in V_M$. Start with
$J_Mv=(v,0)$. If $t$ is not one of the selected indices, i.e., $d(t) =d(t+1)$, define
$J_{t+1}=J_tR_t$, where $R_t:V_{t+1}\to V_t$ is the
canonical retraction given by
$R_t(x,a)=x$. If $t=s_k$, i.e., $d(t+1) = d(t) +1$, define
\[
 J_{t+1}(x,a)=\bigl(J_t(x-a_1p),a_1\bigr),
\]
where the last coordinate is appended to the $\ell_1$-part. In the
second case,
\[
 \norm{J_{t+1}(x,a)}
 \le\norm{x-a_1p}_{V_t}+|a_1|
 \le\norm{(x,a)}_{V_{t+1}};
\]
in the first case contractivity follows from that of $J_t$ and $R_t$.
Both definitions preserve the stated action on $V_M$.

More explicitly, induction gives, for every $x\in V_t$,
\begin{equation}\label{eq:Jselected}
 J_tx=
 \left(\rho_Mx-p\sum_{s_k<t}\lambda_k(x),\,
              (\lambda_k(x))_{s_k<t}\right).
\end{equation}

The case $t = M$ is trivial. Suppose that
$t\notin\{s_1,\ldots,s_r\}$.
By definition, we have
\[
J_{t+1}=J_tR_t, \quad\quad d(t+1)=d(t).
\]
Let $z=(x,a)\in V_{t+1}$.
The canonical retractions are
compatible, so
\[
\rho_Mz=\rho_M(R_tz)=\rho_Mx.
\]
Therefore the formula (7.3) also holds for $J_{t+1}$ if it holds for $J_t$. Suppose that $t=s_j$ and again let
$z=(x,a)\in V_{t+1}$. Then
\begin{align*}
J_{t+1}z
&=\bigl(J_t(x-a_1p),a_1\bigr)\\
&=\left(
\rho_M(x-a_1p)
-p\sum_{k<j}\lambda_k(x-a_1p),
(\lambda_k(x-a_1p))_{k<j},a_1
\right)\\
&=\left(
\rho_Mx-a_1p
-p\sum_{k<j}
\bigl(\lambda_k(x)-a_1\lambda_k(p)\bigr),
(\lambda_k(x)-a_1\lambda_k(p))_{k<j},
a_1
\right)\\
&=\left(
\rho_Mx-a_1p
-p\sum_{k<j}\lambda_k(x),
(\lambda_k(x))_{k<j},a_1
\right)\\
&=\left(
\rho_Mz
-p\left(\sum_{k<j}\lambda_k(z)
+\lambda_j(z)\right),
(\lambda_k(z))_{k\le j}
\right)\\
&=\left(
\rho_Mz-p\sum_{s_k<t+1}\lambda_k(z),
(\lambda_k(z))_{s_k<t+1}
\right),
\end{align*}
since $\lambda_k(p) = 0$. At $t=s_r+1$, contractivity and \eqref{eq:Jselected} imply
\[
 \norm{\rho_Mx-p\sum_{k=1}^r\lambda_k(x)}
       +\sum_{k=1}^r|\lambda_k(x)|\le\norm x
 \qquad(x\in V_{s_r+1}).
\]
Discarding the first term gives \eqref{eq:coordinates} on this stage.
For $x\in X$, apply that inequality to $\rho_{s_r+1}x$.
Compatibility gives
$\lambda_k(\rho_{s_r+1}x)=\lambda_k(x)$, and
$\norm{\rho_{s_r+1}x}\le\norm x$. The final assertion follows by
the triangle inequality for scalar values.
\end{proof}

\begin{proposition}\label{prop:DPY}
The space $Y$ has the Daugavet property.
\end{proposition}

\begin{proof}
First fix $f,g\in S_{Y_N}$ and $\eta>0$. Since $\mathcal P$ is dense
in $S_X$ by \eqref{eq:recurrence} and $\norm g_{X^*}=1$, choose
$p\in\mathcal P$ with $g(p)>1-\eta$. Choose $M\ge N$ with
$p\in V_M$. The restriction of $f$ to $V_M$ still has norm one:
it is at most one, and already has norm one on $V_N$ by the isometry
of $\rho_N^*$.

For any positive integer $r$, the recurrence in
\eqref{eq:recurrence} gives indices
$M\le s_1<\cdots<s_r$ with $p_{s_k}=p$. Let $\lambda_k$ be as in
Lemma~\ref{lem:coordinates}, and put
\[
 t=1-f(p)\in[0,2],\qquad h_k=f+t\lambda_k\in Y.
\]
On $V_{s_k+1}=G(V_{s_k},p,\mu_{s_k})$, the old functional $f$ ignores
the new coordinates. Therefore
\[
 h_k(x,a)=f(x)+(1-f(p))a_1=f(x-a_1p)+a_1,
\]
and hence
\[
 |h_k(x,a)|\le\norm{x-a_1p}+|a_1|
              \le\norm{(x,a)}_{V_{s_k+1}}.
\]
Both $f$ and $\lambda_k$ factor through $\rho_{s_k+1}$, therefore also $h_k$, so this
inequality gives $\norm{h_k}_{X^*}\le1$. On $V_M$, one has
$h_k=f$, because $\lambda_k$ vanishes there. Thus $\norm{h_k}=1$.

Let $q_k=(p,u_1)\in V_{s_k+1}$. Lemma~\ref{lem:extension} gives
$\norm{q_k}=1$. Since $f,g\in Y_N$ and $N\le s_k$,
\[
 h_k(q_k)=f(p)+t=1,\qquad g(q_k)=g(p)>1-\eta.
\]
It follows that
\[
 \norm{g+h_k}\ge(g+h_k)(q_k)>2-\eta.
\]
On the other hand, Lemma~\ref{lem:coordinates} gives
\begin{equation}\label{eq:Ymean}
 \norm{\frac1r\sum_{k=1}^rh_k-f}
 =\frac tr\norm{\sum_{k=1}^r\lambda_k}\le\frac2r.
\end{equation}
Thus these two estimates hold with arbitrarily small errors for
normalized vectors in a common finite dual stage.

For arbitrary $f,g\in S_Y$ and $\eta>0$, choose $0<\tau<\eta/3$ and
$f_0,g_0\in S_{Y_N}$ with
$\norm{f-f_0}<\tau$ and $\norm{g-g_0}<\tau$. This is possible by
\eqref{eq:Y}. Apply the preceding construction to $f_0,g_0$ with
tolerance $\tau$ and an integer $r$ such that $2/r<\tau$. Then
\[
 \norm{\frac1r\sum_{k=1}^rh_k-f}
 \le\norm{\frac1r\sum_{k=1}^rh_k-f_0}+\norm{f_0-f}
 <2\tau<\eta,
\]
while
\[
 \norm{g+h_k}\ge\norm{g_0+h_k}-\norm{g-g_0}
 >2-2\tau>2-\eta.
\]
The proof is therefore finished by Theorem 1.1 (3).
\end{proof}

\section{The projective tensor square}

Put $Z_\pi=Y\tp Y$ and regard $Z_\pi^*$ as the space of bounded
bilinear forms on $Y\times Y$. Let $e_i^*\in H^*$ be the coordinate
functionals and define
\begin{equation}\label{eq:FB}
 u_i=\rho_0^*e_i^*\in Y_0,\qquad
 F=\tfrac12(u_1\otimes u_1+u_2\otimes u_2),\qquad
 B_0(f,g)=\sum_{i=1}^2f(e_i)g(e_i).
\end{equation}
Restriction $f\mapsto f|_H$ is a contraction from $Y$ to $H^*$.
The Cauchy--Schwarz inequality therefore gives $\norm{B_0}\le1$.
Also $\norm{u_i}=1$ and $u_i(e_j)=\delta_{ij}$, so
$B_0(u_i,u_i)=1$ and therefore
\begin{align*}
B_0(F)
&=B_0\left(\frac12(u_1\otimes u_1+u_2\otimes u_2)\right)\\
&=\frac12\left(B_0(u_1,u_1)+B_0(u_2,u_2)\right)\\
&=\frac12\sum_{i=1}^{2}\sum_{j=1}^{2}u_i(e_j)^2\\
&=\frac12\sum_{i=1}^{2}\sum_{j=1}^{2}\delta_{ij}^{2}=\frac12(1+1)=1,
\end{align*}
by the identification of $Z_\pi^*$ with the bounded bilinear forms on $Y \times Y$. The representation of $F$ in
\eqref{eq:FB} gives $\norm F_\pi\le1$, while evaluation by $B_0$
gives the reverse inequality. Thus
\begin{equation}\label{eq:FBnorms}
 \norm{B_0}=\norm F_\pi=B_0(F)=1.
\end{equation}

\begin{proposition}\label{prop:projectivebound}
For every $B\in B_{Z_\pi^*}$,
\begin{equation}\label{eq:projectivebound}
 \norm{B-B_0}\le2\sqrt{1-B(F)}+122\bigl(1-B(F)\bigr).
\end{equation}
In particular, evaluation at $F$ weak-star strongly exposes $B_0$ in
$B_{Z_\pi^*}$.
\end{proposition}

\begin{proof}
Fix $s\ge0$. Via the isometry $\rho_s^*:V_s^*\to Y_s$, the
restriction of $B$ to $Y_s\times Y_s$ becomes the bilinear form
\[
 \widetilde B_s(a^*,b^*)=B(\rho_s^*a^*,\rho_s^*b^*)
 \qquad(a^*,b^*\in V_s^*).
\]
Since $V_s$ is finite-dimensional, there is a unique tensor
$w_s\in V_s\otimes V_s$ with
\[
 (a^*\otimes b^*)(w_s)=\widetilde B_s(a^*,b^*)
 \quad(a^*,b^*\in V_s^*).
\]
Taking the supremum over the two unit balls gives exactly
\begin{equation}\label{eq:bilinearnorm}
 \norm{w_s}_\eps=\norm{B|_{Y_s\times Y_s}}\le1.
\end{equation}
The analogous tensor for $B_0|_{Y_s\times Y_s}$ is $z_0$, since
$(\rho_s^*a^*)(e_i)=a^*(e_i)$.

Write
\[
w_s=\sum_{j=1}^m x_j\otimes y_j,
\qquad x_j,y_j\in V_s.
\]

For the compressed trace, note first that
$\rho_s^*P_s^*e_i^*=\rho_0^*e_i^*=u_i$, because
$P_s\rho_s=P$. By the matrix identification \eqref{eq:matrix},
\begin{align*}
\Phi_s(w_s)
&=\frac12\tr\bigl((P_s\otimes P_s)w_s\bigr)=\frac12\sum_{i=1}^2
e_i^*\bigl(P_sw_sP_s^*e_i^*\bigr)=\frac12\sum_{i=1}^2
(P_s^*e_i^*)\bigl(w_s(P_s^*e_i^*)\bigr)\\
&=\frac12\sum_{i=1}^2
(P_s^*e_i^*)\left(
\sum_{j=1}^m
(P_s^*e_i^*)(x_j)y_j
\right)=\frac12\sum_{i=1}^2\sum_{j=1}^m
(P_s^*e_i^*)(x_j)
(P_s^*e_i^*)(y_j)\\
&=\frac12\sum_{i=1}^2
(P_s^*e_i^*\otimes P_s^*e_i^*)(w_s)=\frac12\sum_{i=1}^2
\widetilde B_s(P_s^*e_i^*,P_s^*e_i^*)\\
&=\frac12\sum_{i=1}^2
B(\rho_s^*P_s^*e_i^*,
  \rho_s^*P_s^*e_i^*)=\frac12\sum_{i=1}^2 B(u_i,u_i)\\
&=B\left(
\frac12\sum_{i=1}^2u_i\otimes u_i
\right)\\
&=B(F).
\end{align*}
The trace defect is thus independent of $s$. Proposition~\ref{prop:finitebound}
and \eqref{eq:bilinearnorm}, applied also to the difference, yield
\begin{equation}\label{eq:restrbound}
 \norm{(B-B_0)|_{Y_s\times Y_s}}
 =\norm{w_s-z_0}_\eps
 \le2\sqrt{1-B(F)}+122\bigl(1-B(F)\bigr).
\end{equation}

For completeness, if $C$ is any bounded bilinear form on $Y\times Y$,
then
\[
 \norm C=\sup_s\norm{C|_{Y_s\times Y_s}}.
\]
One inequality follows from restriction. For the other, fix
$f,g\in S_Y$. By density and normalization choose $f_j,g_j$ in the
unit spheres of a common stage $Y_{s_j}$, with $f_j\to f$ and
$g_j\to g$. Bilinear continuity gives
$C(f,g)=\lim_j C(f_j,g_j)$, whose absolute value is bounded by the
displayed supremum. Now take the supremum over $f,g\in S_Y$.
Applying this observation to \eqref{eq:restrbound} proves
\eqref{eq:projectivebound}. Since the exposing functional is evaluation
at the element $F\in Z_\pi$, the strong exposure is weak-star strong
exposure in the dual ball.
\end{proof}

\begin{corollary}\label{cor:projective}
The operator $T_\pi z=-B_0(z)F$ on $Z_\pi$ has rank one,
$\norm{T_\pi}=1$, and
$\norm{\id_{Z_\pi}+T_\pi}\le199/100$.
\end{corollary}

\begin{proof}
The norm and rank assertions follow from \eqref{eq:FBnorms}. For
$B\in Z_\pi^*$ one has
\[
 T_\pi^*B=-B(F)B_0.
\]
Let $\norm B\le1$. If $|B(F)|\le99/100$, then
\[
 \norm{B+T_\pi^*B}\le1+|B(F)|\le199/100.
\]
Otherwise put $C=\operatorname{sign}(B(F))B$, so that
$C(F)=|B(F)|>99/100$. By Proposition~\ref{prop:projectivebound},
\[
 \norm{C-B_0}\le2\sqrt{1/100}+122/100=71/50.
\]
It follows that
\[
 \norm{B+T_\pi^*B}
 =\norm{C-C(F)B_0}
 \le\norm{C-B_0}+1-C(F)\le143/100.
\]
Consequently
\[
 \norm{\id_{Z_\pi}+T_\pi}
 =\norm{\id_{Z_\pi^*}+T_\pi^*}\le199/100.
\]
\end{proof}

Propositions~\ref{prop:DPX} and~\ref{prop:DPY}, together with
Corollaries~\ref{cor:injective} and~\ref{cor:projective}, prove
Theorem~\ref{thm:main}.

\section{Relation with the operator Daugavet properties}

The \emph{operator Daugavet property} (ODP), introduced in
\cite[Definition~4.1]{RTV}, means that, given
$x_1,\ldots,x_n\in S_E$, $\eta>0$, and a slice $S$ of $B_E$, there is
$h\in S$ such that for every $q\in B_E$ there is an operator $A:E\to E$
with $\norm A\le1+\eta$, $\norm{Ax_i-x_i}<\eta$ for every $i$, and
$Ah=q$. The weak operator Daugavet property of \cite[Definition~5.2]{MRZ}
allows $h$ to depend on $q$ and requires only $\norm{Ah-q}<\eta$.

\begin{proposition}\label{prop:ODP}
The space $X$ has the ODP. The space $Y$ does not have the weak
operator Daugavet property.
\end{proposition}

\begin{proof}
Fix $x_1,\ldots,x_n\in S_X$, $\eta>0$, and a slice
$S=\Slice(B_X,x^*,\alpha)$. Choose $v\in S_X$ with
$x^*(v)>1-\alpha/4$. At a sufficiently late outer stage $E_N$, we
can find $y_i,v_0\in E_N$ with $\norm{x_i-y_i}<\eta/3$,
$v_0\in S_{E_N}$ and $x^*(v_0)>1-\alpha/2$.
Choose $n_0\ge N$ with $\eps_{n_0}+2/m_{n_0}<\alpha/2$, and choose
$p\in N_{n_0}$ with $\norm{p-v_0}<\eps_{n_0}$. At the extension
$V_{s+1}=G(V_s,p,m_{n_0})$ corresponding to this point, the vectors
$h_j=(p,u_j)$ have an average at distance at most
$\eps_{n_0}+2/m_{n_0}$ from $v_0$ by Lemma 2.1. Hence some $h_j$ belongs to $S$.
Fix this $j$.

For any $q\in B_X$, define
\[
 A_q^0:V_{s+1}\to X,\qquad
 A_q^0(x,a)=x-a_jp+a_jq.
\]
The extension norm gives
\[
 \norm{A_q^0(x,a)}\le\norm{x-a_jp}+|a_j|
                    \le\norm{(x,a)}_{V_{s+1}}.
\]
Thus $A_q=A_q^0\rho_{s+1}:X\to X$ is contractive. It fixes every
vector of $V_s$, sends $h_j$ to $q$, and satisfies
\[
 \norm{A_qx_i-x_i}\le2\norm{x_i-y_i}<\eta.
\]
The choice of $h_j$ is independent of $q$, proving the ODP of $X$.

If $Y$ had the weak operator Daugavet property, then
\cite[Theorem~5.4]{MRZ} would imply that $Y\tp Y$ has that property,
and hence the Daugavet property by \cite[Remark~5.3]{MRZ}. This contradicts
Corollary~\ref{cor:projective}.
\end{proof}

There are also natural examples of spaces with the Daugavet property and without the weak operator Daugavet property.
We use the following criterion:
\begin{lemma}[{\cite[Lemma~2.2]{Hamad}}]\label{lem:hamad}
    Let $E$ be a Banach space with the WODP, and let $\delta,\varepsilon>0$.
    Suppose that $x_1,x_2\in S_E$ satisfy
    \[
        B_E\!\left(x_1,\tfrac12\right)
        \cap B_E\!\left(x_2,\tfrac12\right)\neq\varnothing,
    \]
    and that every
    \[
        z\in B_E\!\left(x_1,\tfrac12+\tfrac32\varepsilon\right)
        \cap B_E\!\left(x_2,\tfrac12+\tfrac32\varepsilon\right)
    \]
    satisfies $\|z+x_1\|>2-\delta$.
    Then, for every slice $S$ of $B_E$, there exists $s\in S$ such that
    \[
        (1+\varepsilon)\|s+u\|>2-\delta-\varepsilon
    \]
    for every
    $u\in B_E(x_1,\tfrac12)\cap B_E(x_2,\tfrac12)$.
\end{lemma}
\begin{proposition}\label{prop:circle-wodp}
    Let
    \[
        M=\mathbb{R}/(3\mathbb{Z}),
        \qquad
        d([s],[t])=\min_{k\in\mathbb{Z}}|s-t+3k|,
    \]
    with distinguished point $[0]$. The Banach space
    \[
        X=\left\{
            f\in L_\infty(0,3;\mathbb{R}):
            \int_0^3 f(t)\,dt=0
        \right\},
        \qquad \|f\|=\|f\|_\infty,
    \]
    is linearly isometric to $\operatorname{Lip}_0(M)$ and does not have the
    WODP while having the Daugavet property. Consequently, $\operatorname{Lip}_0(M)$ does not have the WODP either.
    Here $\operatorname{Lip}_0(M)$ is the space of real-valued Lipschitz
    functions on $M$ vanishing at $[0]$, equipped with the Lipschitz norm.

    The isometry is given by
    \[
        D\colon\operatorname{Lip}_0(M)\longrightarrow X,
        \qquad DF=\widetilde F'\big|_{(0,3)},
    \]
    where $\widetilde F(t)=F([t])$ is the $3$-periodic lift of $F$, and the
    derivative is understood almost everywhere.
\end{proposition}

\begin{proof}
The isometry between $X$ and $\operatorname{Lip}_0(M)$ is clear by
\[
 f\longmapsto F_f,\qquad
 F_f([t])=\int_0^t f(r)\,dr,
 \quad 0\le t\le 3,
\]
and $\operatorname{Lip}_0(M)$ has the Daugavet property, since $M$ is geodesic \cite{GLPRZ}.
    To prove that $X$ does not have the WODP, set
    $A=(0,1)$, $B=(1,2)$, and $C=(2,3)$. We use $(a,b,c)$ to denote the
    function taking the corresponding constant values on these intervals.
    In particular,
    \[
        x_1=(1,-1,0),\qquad x_2=(1,0,-1)
    \]
    belong to $S_X$. The intersection
    $B_X(x_1,\tfrac12)\cap B_X(x_2,\tfrac12)$ is nonempty, since
    \[
        \frac{x_1+x_2}{2}=(1,-\tfrac12,-\tfrac12)
        \in B_X(x_1,\tfrac12)\cap B_X(x_2,\tfrac12).
    \]

    Choose
    \[
        \varepsilon=\frac1{20},\qquad
        \delta=\frac15,\qquad
        \frac12+\frac32\varepsilon=\frac{23}{40},
    \]
    and let
    $z\in B_X(x_1,\tfrac{23}{40})\cap B_X(x_2,\tfrac{23}{40})$.
    On $B$, the two ball conditions give
    \[
        -\frac{63}{40}\leq z\leq-\frac{17}{40},
        \qquad
        -\frac{23}{40}\leq z\leq\frac{23}{40}
        \quad\text{almost everywhere}.
    \]
    The same inequalities, with their roles interchanged, hold on $C$.
    Consequently,
    \[
        -\frac{23}{40}\leq z\leq-\frac{17}{40}
        \quad\text{almost everywhere on }B\cup C.
    \]
    As $z\in X$ and $|B|=|C|=1$, it follows that
    \[
        \int_Az=-\int_{B\cup C}z\geq\frac{17}{20}.
    \]
    Since $|A|=1$ and $x_1=1$ on $A$, we obtain
    \begin{equation}\label{eq:circle-hamad-hypothesis}
        \|z+x_1\|_\infty
        \geq\int_A(z+x_1)
        =1+\int_Az
        \geq\frac{37}{20}
        >\frac95=2-\delta.
    \end{equation}
    Thus the geometric hypothesis of Lemma~\ref{lem:hamad} is satisfied.

    Consider the functional $\varphi\in X^*$ defined by
    $\varphi(f)=-\int_Af$. Clearly, $\|\varphi\|\leq1$.
    Since $(-1,\tfrac12,\tfrac12)\in S_X$ and
    $\varphi(-1,\tfrac12,\tfrac12)=1$, we have $\|\varphi\|=1$.
    Therefore
    \[
        S=\{s\in B_X:\varphi(s)>\tfrac12\}
         =\left\{s\in B_X:\int_As<-\frac12\right\}
    \]
    is a nonempty slice of $B_X$.

    Fix $s\in S$, put
    \[
        a=\int_As\in[-1,-\tfrac12),
        \qquad
        c=\frac1{2(1-a)}\in[\tfrac14,\tfrac13),
    \]
    and define
    \[
        u(t)=
        \begin{cases}
            \displaystyle\frac12+c\bigl(1-s(t)\bigr),&t\in A,\\[1mm]
            -\displaystyle\frac12,&t\in B\cup C.
        \end{cases}
    \]
    We have
    \[
        \int_Au=\frac12+c(1-a)=1,
        \qquad
        \int_{B\cup C}u=-1.
    \]
    Hence $u\in X$. Moreover, since $-1\leq s\leq1$ almost everywhere,
    \[
        \frac12\leq u\leq\frac12+2c
        <\frac76<\frac32
        \quad\text{almost everywhere on }A.
    \]
    On $A$ both $x_1$ and $x_2$ equal $1$, while on $B\cup C$ the function
    $u$ equals $-\tfrac12$ and each of $x_1,x_2$ takes values in $\{-1,0\}$.
    Thus
    \[
        \|u-x_1\|_\infty\leq\frac12,
        \qquad
        \|u-x_2\|_\infty\leq\frac12,
    \]
    so $u\in B_X(x_1,\tfrac12)\cap B_X(x_2,\tfrac12)$.

    On $A$ we have
    \[
        s+u=\frac12+c+(1-c)s.
    \]
    Since $1-c>0$ and $-1\leq s\leq1$ almost everywhere, this gives
    \[
        -\frac12+2c\leq s+u\leq\frac32.
    \]
    The lower bound is nonnegative because $c\geq\tfrac14$.
    On $B\cup C$, we have $|s+u|=|s-\tfrac12|\leq\tfrac32$.
    Consequently,
    \begin{equation}\label{eq:circle-small-sum}
        \|s+u\|_\infty\leq\frac32.
    \end{equation}

    Suppose now that $X$ has the WODP. By
    \eqref{eq:circle-hamad-hypothesis}, Lemma~\ref{lem:hamad} provides an
    $s\in S$ such that
    \[
        (1+\varepsilon)\|s+u\|_\infty
        >2-\delta-\varepsilon=\frac74
    \]
    for every $u\in B_X(x_1,\tfrac12)\cap B_X(x_2,\tfrac12)$.
    For this $s$, however, the preceding construction gives an element
    $u$ in the same intersection satisfying \eqref{eq:circle-small-sum}.
    Therefore
    \[
        (1+\varepsilon)\|s+u\|_\infty
        \leq\frac{21}{20}\cdot\frac32
        =\frac{63}{40}
        <\frac74,
    \]
    a contradiction. 
\end{proof}

\section{Isomorphic result}
For the main result of this section, we need the following theorem from \cite{Shvidkoy}.

\begin{theorem}
Suppose that $E$ has the Daugavet property and $F$ is a closed subspace of $E$.
If the quotient $E/F$ has the Radon--Nikod\'ym property, then $F$ also has the
Daugavet property with its inherited norm.
    
\end{theorem}
\begin{theorem}
    Let $X$ be a nonzero Banach space and put
    \[
        \kappa=\max\{\operatorname{dens}(X),\aleph_0\},
    \]
    where $\operatorname{dens}(X)$ denotes the density character of $X$.
    Let $\Gamma$ be a set such that $|\Gamma|>\kappa$.
    Then the Banach space $Y=X\oplus_2\ell_2(\Gamma)$ has the following
    properties:
    \begin{enumerate}
        \item $X$ embeds linearly and isometrically into $Y$.
        \item No equivalent norm on $Y$ has the Daugavet property.
    \end{enumerate}
    Consequently, $Y$ is not isomorphic to any Banach space with the
    Daugavet property.
\end{theorem}

\begin{proof}
    Equip $Y=X\oplus_2\ell_2(\Gamma)$ with the norm
    \[
        \|(x,h)\|_0=\left(\|x\|_X^2+\|h\|_2^2\right)^{1/2}.
    \]
    The canonical map
    \[
        J\colon X\longrightarrow Y,
        \qquad J(x)=(x,0),
    \]
    is a linear isometry. We identify $X$ with $J(X)$ throughout the proof.

    Suppose, towards a contradiction, that $p$ is an equivalent norm on
    $Y$ for which $(Y,p)$ has the Daugavet property. Since equivalent norms
    induce the same topology, the unit sphere of $(X,p|_X)$ has a dense
    subset $D$ with $|D|\leq\kappa$. For each $d\in D$, the Hahn--Banach
    theorem gives $F_d\in(Y,p)^*$ such that
    \[
        \|F_d\|_{p,*}=1,
        \qquad F_d(d)=1.
    \]
    This family is norming for $(X,p|_X)$; that is,
    \begin{equation}\label{eq:daugavet-embedding-norming}
        p(x)=\sup_{d\in D}|F_d(x)|
        \qquad(x\in X).
    \end{equation}
    Indeed, $|F_d(x)|\leq p(x)$ for every $d\in D$ and $x\in X$.
    Conversely, if $p(x)=1$ and $\varepsilon>0$, choose $d\in D$ with
    $p(x-d)<\varepsilon$. Then
    \[
        |F_d(x)|
        \geq |F_d(d)|-|F_d(x-d)|
        \geq 1-p(x-d)
        >1-\varepsilon.
    \]
    Letting $\varepsilon\to0$ gives
    \eqref{eq:daugavet-embedding-norming} for unit vectors, and the general
    case follows by homogeneity.

    Identify $\ell_2(\Gamma)$ with the closed subspace
    $H=\{0\}\oplus\ell_2(\Gamma)$ of $Y$. Since $p$ is equivalent to
    $\|\cdot\|_0$, each restriction $F_d|_H$ is continuous for the original
    Hilbert space norm. By the Riesz representation theorem, it is
    represented by a vector $a_d\in\ell_2(\Gamma)$. Every such vector has
    at most countable support, so
    \[
        S=\bigcup_{d\in D}\operatorname{supp}(a_d)
        \qquad\text{satisfies}\qquad
        |S|\leq |D|\cdot\aleph_0\leq\kappa.
    \]
    As $|\Gamma|>\kappa$, we can choose $\gamma_0\in\Gamma\setminus S$.
    Set $e=(0,e_{\gamma_0})$, where $e_{\gamma_0}$ is the corresponding
    canonical unit vector of $\ell_2(\Gamma)$. By construction,
    \[
        F_d(e)=0
        \qquad(d\in D).
    \]
    For every $x\in X$ and $t\in\mathbb{R}$, it follows from
    \eqref{eq:daugavet-embedding-norming} that
    \[
        \begin{aligned}
            p(x)
                &=\sup_{d\in D}|F_d(x)| \\
                &=\sup_{d\in D}|F_d(x+te)| \\
                &\leq p(x+te).
        \end{aligned}
    \]
    Thus, on the closed subspace $Z=X\oplus\mathbb{R}e$, the projection
    \[
        Q\colon Z\longrightarrow Z,
        \qquad Q(x+te)=x,
    \]
    is contractive. Since $Q$ is the identity on the nonzero subspace $X$,
    we have $\|Q\|_p=1$.

    The quotient $Y/Z$, equipped with the quotient norm induced by
    $\|\cdot\|_0$, is isometric to
    $\ell_2(\Gamma\setminus\{\gamma_0\})$. Equivalent norms on $Y$ induce
    equivalent quotient norms, so $(Y,p)/Z$ is isomorphic to this Hilbert
    space and is therefore reflexive. In particular, it has the
    Radon--Nikod\'ym property. The hereditary theorem stated above now
    implies that $(Z,p|_Z)$ has the Daugavet property.

    Consider $P=I_Z-Q$. Since
    \[
        P(x+te)=te
        \qquad(x\in X,\ t\in\mathbb{R}),
    \]
    the operator $P$ is a nonzero rank-one projection. Moreover,
    $\|P\|_p\geq1$, because $P(e)=e$. Applying the Daugavet equation in
    $(Z,p|_Z)$ to the rank-one operator $-P$ gives
    \[
        1=\|Q\|_p=\|I_Z-P\|_p=1+\|P\|_p\geq2,
    \]
    a contradiction. Hence no equivalent norm on $Y$ has the Daugavet
    property.
\end{proof}

\section*{Acknowledgments}
The author used OpenAI's GPT-6 Astra to assist with mathematical exploration,
the development and checking of proofs, and the presentation of arguments,
particularly those concerning injective and projective tensor products.
The author takes full responsibility for the mathematical content and any
remaining errors.


\begin{thebibliography}{99}

\bibitem{Shvidkoy}
    R.~V.~Shvidkoy,
    \emph{Geometric aspects of the Daugavet property},
    Journal of Functional Analysis \textbf{176} (2000), no.~2, 198--212.

\bibitem{GLPRZ}
L.~García-Lirola, A.~Procházka and A.~Rueda Zoca,
\emph{A characterisation of the Daugavet property in spaces of Lipschitz
functions},
Journal of Mathematical Analysis and Applications \textbf{464} (2018),
473--492.
\href{https://doi.org/10.1016/j.jmaa.2018.04.017}
{doi:10.1016/j.jmaa.2018.04.017}.

\bibitem{Hamad}
Samir Hamad,
\emph{On the weak operator Daugavet property and poor subspaces},
preprint, arXiv:2608.11098 [math.FA], 2026.
\url{https://arxiv.org/abs/2608.11098}.

\bibitem{KKW}
Vladimir Kadets, Nigel Kalton, and Dirk Werner,
\emph{Remarks on rich subspaces of Banach spaces},
Studia Math. \textbf{159} (2003), no.~2, 195--206.
\href{https://doi.org/10.4064/sm159-2-3}{doi:10.4064/sm159-2-3}.

\bibitem{KSSW}
Vladimir M. Kadets, Roman V. Shvidkoy, Gleb G. Sirotkin, and Dirk Werner,
\emph{Banach spaces with the Daugavet property},
Trans. Amer. Math. Soc. \textbf{352} (2000), no.~2, 855--873.
\href{https://doi.org/10.1090/S0002-9947-99-02377-6}{doi:10.1090/S0002-9947-99-02377-6}.

\bibitem{MRZ}
Miguel Mart\'in and Abraham Rueda Zoca,
\emph{Daugavet property in projective symmetric tensor products of Banach spaces},
Banach J. Math. Anal. \textbf{16} (2022), article~35.
\href{https://doi.org/10.1007/s43037-022-00186-6}{doi:10.1007/s43037-022-00186-6}.

\bibitem{RZ}
Abraham Rueda Zoca,
\emph{Weak operator Daugavet property and weakly open sets in tensor product spaces},
Rev. Mat. Complut. \textbf{39} (2026), 149--164.
\href{https://doi.org/10.1007/s13163-025-00531-x}{doi:10.1007/s13163-025-00531-x}.

\bibitem{RTV}
Abraham Rueda Zoca, Pedro Tradacete, and Ignacio Villanueva,
\emph{Daugavet property in tensor product spaces},
J. Inst. Math. Jussieu \textbf{20} (2021), no.~4, 1409--1428.
\href{https://doi.org/10.1017/S147474801900063X}{doi:10.1017/S147474801900063X}.

\bibitem{Ryan}
Raymond A. Ryan,
\emph{Introduction to Tensor Products of Banach Spaces},
Springer Monographs in Mathematics, Springer-Verlag, London, 2002.
\href{https://doi.org/10.1007/978-1-4471-3903-4}{doi:10.1007/978-1-4471-3903-4}.

\bibitem{Werner}
Dirk Werner,
\emph{Recent progress on the Daugavet property},
Irish Math. Soc. Bull. \textbf{46} (2001), 77--97.
\href{https://doi.org/10.33232/BIMS.0046.77.97}{doi:10.33232/BIMS.0046.77.97}.
\bibitem{KadetsMartinRuedaZocaWerner}
V. Kadets, M. Mart\'in, A. Rueda Zoca,
and D. Werner,
\emph{Banach Spaces with the Daugavet Property},
preliminary version, May 2025.
\end{thebibliography}
\end{document}